\documentclass[12pt]{article}
\usepackage{times}                               
\usepackage{amsmath}

\title
{
On a Combinatorial Identity of Djakov and Mityagin
}
\author{Steven J. Rosenberg\\
\small Department of Mathematics and Computer Science\\
\small University of Wisconsin -- Superior,\\
\small Superior, Wisconsin, USA\\
\small \texttt{srosenbe@uwsuper.edu}
}
\date{
}

\begin{document}
\maketitle

\begin{abstract}
Consider a pyramid made out of unit cubes arranged in square horizontal layers,
with a ledge of one cube's length around the perimeter of each layer.
For any natural number $k$, we can count the number of ways of choosing
$k$ unit cubes from the pyramid such that no two cubes are in the same horizontal layer;
we can also count the number of ways of choosing $k$ unit cubes from the pyramid
such that no two cubes come from the same vertical slice
(taken parallel to a fixed edge of the pyramid) {\it or} from two adjacent slices.
Djakov and Mityagin first established, using functional analysis,
that these two quantities are always equal
(the enumerative interpretation given here is due to
Thomas Kalinowski).
Don Zagier supplied the first combinatorial proof of this result.
We provide a new, more natural combinatorial proof.
\end{abstract}

The authors of \cite{borwein}
conjectured (\cite[Conjecture 1]{borweinArxiv1};
\cite[Theorem 4]{borweinArxiv2, borwein})
that the following identity holds
for all positive integers $n$ and integers $k$:
\begin{equation}\label{conj:Borwein}
\sum_{J\in\binom{[n]}{k}^{*}} \prod_{j\in J} j^2 = \sum_{J\in\binom{[n]}{k}^{**}} \prod_{j\in J} j(n + 1 - j)
\end{equation}
where $[n] := \{1, 2, 3, \dots, n\}$ as usual,
$\binom{[n]}{k}^{*}$ denotes the collection of $k$-subsets of $[n]$ all of whose elements
are congruent to $n$ modulo $2$,
and $\binom{[n]}{k}^{**}$ denotes the collection of $k$-subsets of $[n]$ which do not contain
two consecutive integers.
Don Zagier gave a combinatorial proof of this result in the appendix to \cite{borweinArxiv2, borwein}.
It turned out that a few years prior to this, the same identity had
been established, with
a less direct, non-combinatorial proof, in the context of
functional analysis
(\cite{mit1, mit2}).
In the present note, we establish the identity in
the most straightforward fashion yet,
by calculating the characteristic
polynomial of a certain matrix in two different ways.

Given two sequences $(a_j)$ and $(b_j)$
of real numbers, consider the
$(n+1)$ by $(n+1)$ matrix
\begin{equation}
M_n = \left[
\begin{array}{lllllll}
0&b_1&0&0&\cdots&0&0\\
a_1&0&b_2&0&\cdots&0&0\\
0&a_2&0&b_3&\cdots&0&0\\
0&0&a_3&0&\cdots&0&0\\
\vdots&\vdots&\vdots&\vdots&\ddots&\vdots&\vdots\\
0&0&0&0&\cdots&0&b_n\\
0&0&0&0&\cdots&a_n&0
\end{array}
\right],
\end{equation}
where the entries in the diagonal below the main diagonal are
$a_1$, $a_2$, $\dots$, $a_n$,
the entries in the diagonal above the main diagonal are
$b_1$, $b_2$, $\dots$, $b_n$,
and the remaining entries are $0$.
Let $\chi_n(x) = \det(xI - M)$ be the characteristic
polynomial of $M_n$ in the variable $x$.
We claim that for all $n\geq 1$ and all $k\in{\bf Z}$,
the coefficient $d_{k,n}$ of $x^{n+1-2k}$ in $\chi_n(x)$
satisfies
\begin{equation}\label{eqn:RHS}
d_{k,n} = (-1)^k\sum_{J\in\binom{[n]}{k}^{**}} \prod_{j\in J} a_jb_j.
\end{equation}
One way to see this is by observing that both sides of Equation \ref{eqn:RHS}
satisfy the recurrence
\begin{equation}
r_{k,n} = r_{k,n-1} - a_{n}b_{n}r_{k-1,n-2}
\end{equation}
for $n\geq 3$, and they agree for $n=1,2$.
For the left-hand side of Equation \ref{eqn:RHS},
the recursion follows by expanding the determinant
about the last row.
For the right-hand side, it follows by considering
the terms of total degree $k$ with respect to the variables
$x_1$, $x_2$, $\dots$, $x_n$ in the product
\begin{equation}
\prod_{j=1}^n (1 - x_ja_jb_j),
\end{equation}
after eliminating terms in which two consecutive indices appear.

Specializing to the case $a_i = i$ and $b_i = n + 1 - i$, the matrix
$M_n$ becomes what Taussky and Todd term the ``Kac matrix" $S_n$
(see \cite{tt} for the long and impressive history of this matrix\footnote{The author wishes to thank R. Brualdi for informing him of this earlier work.};
 though named for Victor Kac, the notation honors Sylvester, who found
 its characteristic polynomial in the 1850s).
We only need the eigenvalues of $S_n$; these are computed, e.g., in
\cite{tt}, but we present an original derivation of the spectrum in what follows.

A non-zero vector $v = \left[v_1\ v_2\ \dots\ v_{n+1}\right]^\text{t}$
is an eigenvector of $S_n$ with eigenvalue $\lambda$ if and only if
the components $v_i$ of $v$ satisfy the relations
\begin{equation}\label{eqn:vRelation}
iv_i + (n-i)v_{i+2} = \lambda v_{i+1}
\end{equation}
for all $i\in\{0, 1, 2, \dots, n\}$;
notice that the coefficients of $v_0$ and $v_{n+2}$ will be $0$,
so this really does work.

Let $d$ be an integer between $0$ and $n$ (inclusive).
Set $v_1 = 1$ (arbitrarily), and
use Equation \ref{eqn:vRelation} to obtain $v_2$, $\dots$, $v_{d+1}$
by setting
\begin{equation}
v_{i+2} = \frac{(n - 2d)v_{i+1} - iv_i}{n - i}
\end{equation}
for $i = 0, 1, \dots, d-1$.
Let $p(x)$ be the unique polynomial of degree at most $d$ with values
$p(i) = v_i$ for $i\in\{1, 2, \dots, d + 1\}$.
Then the equation
\begin{equation}\label{eqn:pRelation}
xp(x) + (n-x)p(x+2) = \lambda p(x+1)
\end{equation}
holds, with $\lambda = n - 2d$, for at least
the $d$ distinct values $0, 1, 2, \dots, d-1$ of $x$.
Now both sides of Equation \ref{eqn:pRelation}
are polynomials of degree at most $d$.
If $\deg(p) < d$, then both sides have degree less than $d$,
so Equation \ref{eqn:pRelation} is an identity for $p$.
If $\deg(p) = d$, then both sides of Equation \ref{eqn:pRelation}
not only have degree $d$, but also have the same leading coefficient
of $(n - 2d)c_d$, where $c_d$ is the leading coefficient of $p$;
subtracting this term from both sides, we find two polynomials of degree
at most $d-1$ which agree at $d$ distinct values of $x$,
so again they are identically equal.
It follows that $\left[p(1)\ p(2)\ \cdots\ p(n+1)\right]^\text{t}$
is an eigenvector of $S_n$ with eigenvalue $n - 2d$.

Since $S_n$ is an $n+1$ by $n+1$ matrix, it follows that
the characteristic polynomial of $S_n$ is
\begin{equation}\label{eqn:LHS}
\chi_n(x) = \prod_{d=0}^n (x - (n - 2d))
= x^\varepsilon\cdot\prod_{\genfrac{}{}{0 pt}{}{1\leq j \leq n}{j \equiv n \pmod{2}}} (x^2 - j^2),
\end{equation}
where $\varepsilon = $
$0$ if $n \equiv 1\pmod{2}$ and
$\varepsilon = 1$ if $n \equiv 0\pmod{2}$.

From Equation \ref{eqn:LHS},
we see that the coefficient of $x^{n + 1 - 2k}$
in $\chi_n(x)$ is
\begin{equation}\label{eqn:final}
d_{k,n} = (-1)^k\sum_{J\in\binom{[n]}{k}^{*}} \prod_{j\in J} j^2 .
\end{equation}
Now Equation \ref{conj:Borwein}
follows by comparing Equations \ref{eqn:RHS} and \ref{eqn:final}.

\end{document}